\documentclass[reqno, letterpaper, oneside]{article}
\usepackage{amsmath,amsfonts,amssymb,amsthm}

\usepackage{bbm}
\usepackage{setspace}
\usepackage{geometry}
\usepackage{enumerate}
\usepackage{verbatim}
\usepackage[colorlinks=false,,linktoc=allbookmarksopen=true,linktocpage,pdftex]{hyperref}
\usepackage[pdftex]{bookmark}
\hypersetup{
	colorlinks=true,       
	linkcolor=blue,        
	citecolor=blue,        
	filecolor=magenta,     
	urlcolor=blue         
}
\usepackage{tabularx} 

\usepackage{mathtools}
\usepackage{cleveref}
\usepackage{enumitem}

\theoremstyle{plain}
\newtheorem{theorem}{Theorem}[section]
\newtheorem{lemma}[theorem]{Lemma}
\newtheorem{proposition}[theorem]{Proposition}
\newtheorem{corollary}[theorem]{Corollary}

\theoremstyle{definition}

\newtheorem{definition}[theorem]{Definition}
\newtheorem{problem}[theorem]{Problem}
\newtheorem{remark}[theorem]{Remark}

\newcommand{\refT}[1]{Theorem~\ref{#1}}
\newcommand{\refC}[1]{Corollary~\ref{#1}}

\newcommand{\refL}[1]{Lemma~\ref{#1}}

\newcommand{\refP}[1]{Proposition~\ref{#1}}

\newcommand{\refD}[1]{Definition~\ref{#1}}

\newcommand{\cA}{\mathcal{A}}

\newcommand{\cC}{\mathcal{C}}
\newcommand{\cD}{\mathcal{D}}
\newcommand{\cE}{\mathcal{E}}
\newcommand{\cF}{\mathcal{F}}

\newcommand{\cH}{\mathcal{H}}

\newcommand{\vol}[1]{\textbf{vol}{(#1)}}
\newcommand{\vold}[2]{\textbf{vol}_{#1}{(#2)}}

\newcommand\lrpar[1]{\left(#1\right)}

\newcommand{\ip}[2]{\left \langle #1,#2 \right \rangle}
\DeclarePairedDelimiter\norm{\lVert}{\rVert}

\newcommand{\lobound}{\partial^0}
\newcommand{\shake}{\phi}
\newcommand\ceil[1]{\lceil#1\rceil}

\newcommand\floor[1]{\lfloor#1\rfloor}

\newcommand\xor{\oplus}
\newcommand\cover[1]{\partial^0#1}
\newcommand{\Proj}{\textbf{Proj}}

\newcommand{\cube}{\mathcal{C}}

\newcommand{\Z}{\mathbb{Z}}
\newcommand{\bP}{\mathbb{P}}
\newcommand{\bE}{\mathbb{E}}
\newcommand{\bS}{\mathbb{S}}
\newcommand{\R}{\mathbb{R}}

\newcommand{\Var}{\mathrm{Var}}

\newcommand{\eps}{\varepsilon}

\let\OLDthebibliography\thebibliography
\renewcommand\thebibliography[1]{
  \OLDthebibliography{#1}
  \setlength{\parskip}{0pt}
  \setlength{\itemsep}{0pt plus 0.3ex}
}

\usepackage{indentfirst}

\title{On the $\ell_0$ Isoperimetric Coefficient for Measurable Sets}
\author{Manuel Fernandez V%
\thanks{School of Mathematics, Georgia Institute of Technology, Atlanta GA~30332, USA. E-mail: {\tt mfernandez39@gatech.edu}.}
}
\date{\today}
\begin{document}
\maketitle  
\begin{abstract}
In this paper we prove that the $\ell_0$ isoperimetric coefficient for axis-aligned cubes, $\psi_{\cube}$, is  $\Theta(n^{-1/2})$ and that the isoperimetric coefficient for any measurable body $K$, $\psi_K$, is of order $O(n^{-1/2})$. As a corollary we deduce that axis-aligned cubes essentially ``maximize'' the $\ell_0$ isoperimetric coefficient: There exists a positive constant $q > 0$ such that
$\psi_K \leq q \cdot \psi_{\cube}$, whenever $\cC$ is an axis-aligned cube and $K$ is any measurable set.
Lastly, we give immediate applications of our results to the mixing time of Coordinate-Hit-and-Run for sampling points uniformly from convex bodies.
\end{abstract}
\section{Introduction}
\par
Efficient sampling from a distribution over a high-dimensional ambient space, say $\R^n$, is a well-studied problem in computer science. Appearing as an intermediary step in many computational algorithms, faster sampling methods automatically speed-up the runtime of algorithms designed for tasks such as volume estimation, optimization, and integration. As a general rule of thumb, one wants to generate a random sample according to a distribution that is ``close'' to the desired distribution in time that is polynomial in the dimension of the ambient space. 
\par 
The most popular approach for sampling methods are those based on Markov chains. Specifically, one designs an ergodic, time-reversible Markov chain whose state space is what you sample from and whose stationary distribution is the distribution you wish to sample with respect to. By the Ergodic theorem, such Markov chains have the property that the starting distribution, no matter its initialization, will converge in distribution  to the chain's stationary distribution. The question then becomes which Markov chains ``mix-rapidly'' (i.e. converge to the stationary distribution in few transitions), with minimal assumptions on the starting distribution.  
\par
For the problem of uniformly sampling points from a convex body $K$ in $\R^n$ two well-known Markov chains are those of the Ball Walk \cite{kannan1997random} and Hit-and-Run \cite{lovasz2004hit}. The state transitions for both Markov chains correspond to picking a uniformly random direction over $\bS^{n-1}$ and, from the current point, moving to a new point in the convex body, along the given direction, according to some criteria. Both Markov chains are well studied and are known to converge to the uniform distribution in time polynomial in $n$.
\par
Another Markov chain, popular in practice, is Gibbs Sampling \cite{turchin1971computation}. Also known as Coordinate Hit-and-Run (CHAR), it is a Markov chain whose state transitions correspond to picking a random coordinate direction and resampling the coordinate of the current point according to the 1-dimensional distribution induced by fixing all other coordinates. Despite its use in practice, strong theoretical guarantees regarding the mixing time of CHAR for uniform sampling from convex bodies were not known until recently.  In work by Laddha and Vempala \cite{laddha2020convergence} the authors showed that CHAR mixes rapidly for a convex body under a ``warm-start'' when the convex body's shape satisfies a certain regularity condition:

\begin{quote}
Let $K$ be a convex body in $\R^n$ containing a unit ball. Let $R^2$ be the expected squared distance of a uniformly random point $K$ from the centroid of $K$, and let $\pi$ denote the uniform distribution on $K$. Let $1_K$ denote the indicator of $K$ and let $\nu_0$ denote the density, with respect to $\pi$, of the initial distribution. If $\norm{\nu_0}_\infty < M$ then the number of steps $T$ needed for CHAR to achieve a density $\nu_T$ satisfying $\norm{\nu_T - 1_K}_{L^2(\pi)} < \eps$ obeys
\begin{equation}\label{thm:warm_start}
T \leq O\lrpar{\frac{M^2R^2n^9}{\eps^2}\log\lrpar{\frac{2M}{\eps}}}.
\end{equation}
\end{quote}
Concurrent work by Narayanan and Srivastava \cite{narayanansrivastava2022} obtained a similar bounded on the mixing time of CHAR from a warm start, albeit with a different regularity condition. Subsequent work by Narayanan, Rajaraman, and Srivastava \cite{narayanan2023sampling} showed that CHAR mixes rapidly even under a ``cold-start'': 
\begin{quote}
Let $K$ be a convex body such that $r \cdot B_\infty \subseteq K \subseteq R \cdot B_\infty$. Let $\pi$ denote the uniform measure on $K$. Let $1_K$ denote the indicator of $K$ and let $\nu_0$ denote the density, with respect to $\pi$, of the initial distribution. If $\norm{\nu_0}_\infty < M$ then the number of steps $T$ needed for CHAR to achieve a density $\nu_T$ satisfying $\norm{\nu_T - 1_K}_{L^2(\pi)} < \eps$ obeys
\begin{equation}\label{thm:cold_start}
T \leq O\lrpar{\frac{n^9}{(r/R)^2}\log\lrpar{\frac{M}{\eps}}}.
\end{equation}
\end{quote}
Both works follow the proof of rapid mixing based on conductance of Markov chains in the continuous setting \cite{lovasz1993random}. 
\par
In order to lower bound the conductance or $s$-conductance of CHAR, both works appealed to a one-step coupling argument, first introduced in \cite{laddha2020convergence}, that reduces the problem to lower bounding the $\ell_0$ isoperimetry of axis-disjoint subsets \ref{def:$l_0$}. This quantity captures the size of $\ell_0$ boundaries \ref{def:cover} of axis-disjoint subsets relative to the size of their complement.
In \cite{laddha2020convergence} the authors gave a lower bound on the $\ell_0$ isoperimetry based on a tiling of the convex body by axis-aligned cubes of a fixed granularity:
\begin{quote}
Let $K$ be a convex body in $\R^n$ containing a unit ball with $R^2 = \bE_K(\norm{x - z_K}^2)$ where $z_k$ is the centroid of $K$. Let $S_1,S_2 \subseteq K$ be two measurable subsets of $K$ such that $S_1$ and $S_2$ are axis-disjoint. Then for any $\eps \geq 0$, the set $S_3 = K\setminus \{S_1 \cup S_2\}$ satisfies
\begin{equation}\label{thm:l_0_iso1}
\vol{S_3} \geq \frac{\eps}{800\cdot n^{3.5}R} \cdot (\min\{\vol{S_1},\vol{S_2}\} - \eps \cdot \vol{K}).
\end{equation}
\end{quote}
Due to the proof method and the use of cubes of fixed granularity, their lower bound has a $-\eps \cdot \vol{K}$ term. This prevents the lower bound from being linear in $\min\{\vol{S_1},\vol{S_2}\}$ and can only be used to prove lower bounds on the $s$-conductance of CHAR. In \cite{narayanan2023sampling} the authors managed to give a lower bound on the $\ell_0$ isoperimetry that was linear in $\min\{\vol{S_1},\vol{S_2}\}$, which itself involved a tiling of the convex body by cubes of varying granularity:
\begin{quote}
Let $K$ be a convex body such that $r \cdot B_\infty \subseteq K \subseteq R \cdot B_\infty$. Let $S_1,S_2 \subseteq K$ be two measurable subsets of $K$ such that $S_1$ and $S_2$ are axis-disjoint. Then the set $S_3 = K\setminus \{S_1 \cup S_2\}$ satisfies
\begin{equation}\label{thm:l_0_iso2}
 \vol{S_3} \geq \Omega \lrpar{\frac{r/R}{n^{3.5}} \min\{\vol{S_1},\vol{S_2}\}}.
\end{equation}
\end{quote}
\par
In both cases the tiling method allowed them to lower bound the $\ell_0$ isoperimetry of convex sets in terms of the $\ell_0$ isoperimetry of axis-aligned cubes.
The following lower bound was given in \cite{laddha2020convergence}, whose proof is an application of the Loomis-Whitney Theorem \cite{loomis1949inequality}.
\begin{quote}
Let $\cube$ denote an axis-aligned cube. $S_1,S_2 \subseteq \cube$ be two measurable subsets of $\cube$ such that $S_1$ and $S_2$ are axis-disjoint.
Then the set $S_3 = \cube \setminus \{S_1 \cup S_2\}$ satisfies
\begin{equation}\label{thm:l_0_cube}
\vol{S_3} \geq \frac{\log 2}{n} \cdot \min\{\vol{S_1},\vol{S_2}\}.
\end{equation}
\end{quote}
An open question left by the authors in \cite{laddha2020convergence} was to determine the correct order of magnitude of the right hand side of the inequality in \ref{thm:l_0_cube}, with any asymptotic improvement implying an automatic improvement to \ref{thm:l_0_iso1} and \ref{thm:l_0_iso2} and therefore the mixing time of CHAR for convex bodies. 
In this paper we answer this question by showing that the $\ell_0$ isoperimetry of axis-aligned cubes is of order $\Theta(n^{-1/2})$.
To state our result precisely we introduce a definition.
\begin{definition}\label{def:iso_coeff}
Let $K \subseteq \R^n$ be a measurable subset of $\R^n$ satisfying $\vol{K} > 0$. We define the $\ell_0$-isoperimetric coefficient of $K$ as
\begin{equation}\label{eq:iso-coef}
\psi_K := \inf_{\underset{S_1,S_2 -\text{axis-disjoint}}{S_1,S_2 \subseteq K}} \frac{\vol{K \setminus (S_1 \cup S_2)}}{\min\{\vol{S_1},\vol{S_2}\}}.
\end{equation}
\end{definition}
\begin{theorem}\label{thm:cube}
$\psi_\cube = \Theta(n^{-1/2})$.
\end{theorem}
Our proof of \refT{thm:cube} uses the $\ell_0$ isoperimetric theorem for the continuous cube, proved by Harper in \cite{harper1999isoperimetric}.  This theorem characterizes the subsets of the continuous unit cube that are  minimizers of $\ell_0$ boundary, see \refD{def:cover}, when we restrict to a class of subsets of a certain structure. We remark that the work of \cite{harper1999isoperimetric} actually addresses the minimizers of vertex-boundary for subsets of the graph product of cliques.
In order to apply Harper's theorem we argue that given any axis-disjoints subsets $S_1,S_2$ of the continuous unit cube, one can transform $S_1,S_2$ into new sets $S_1',S_2'$ such that these sets are also axis-disjoint, the volumes of $S_1',S_2'$ being approximately the same as those of $S_1,S_2$, and with $S_1'$ being a set for which the isoperimetric theorem of Harper applies. The proof of \refT{thm:cube} is then finished by  providing a quantitative lowerbound on the volume of the $\ell_0$ boundaries of the extremizers.

As a corollary to \refT{thm:cube}, the proofs of \eqref{thm:l_0_iso1} and \eqref{thm:l_0_iso2}  from \cite{laddha2020convergence} and \cite{narayanan2023sampling} imply that the right hand side of the inequalities can be increased by a factor of $n^{1/2}$. Consequently, the proofs of \eqref{thm:warm_start} and \eqref{thm:cold_start} from \cite{laddha2020convergence} and \cite{narayanan2023sampling} imply that the right hand side of the inequalities can be reduced by a factor of $n$.
 \begin{corollary}
 Let $T_{\text{warm}},T_\text{cold}$ denote the mixing time of CHAR under the assumptions of \ref{thm:l_0_iso1} and \ref{thm:l_0_iso2} respectively. Then they satisfy
 \begin{equation}\label{eq:improved}
T_{\text{warm}} \leq O\lrpar{\frac{M^2R^2n^8}{\eps^2}\log\lrpar{\frac{2M}{\eps}}}, T_{\text{cold}} \leq O\lrpar{\frac{n^8}{(r/R)^2}\log\lrpar{\frac{M}{\eps}}}.
 \end{equation}
 \end{corollary}
 Given that \refT{thm:cube} gives us the correct order of magnitude of the $\ell_0$-isoperimetric coefficient for axis-aligned cubes it is natural to ask how the isoperimetric coefficient for axis-aligned cubes compares to the isoperimetric coefficient of other types of sets. Towards that end we manage to prove, somewhat surprisingly, a universal upper bound on the $\ell_0$ isoperimetric coefficient for any measurable set.
\begin{theorem}\label{thm:meas}
If $K \subseteq \R^n$ is a measurable set of finite volume then
$\psi_K = O(n^{-1/2})$.
\end{theorem}
 \refT{thm:meas} and \refT{thm:cube} together imply the following: There exists an absolute constant $q > 0$ such that $\psi_K \leq q \cdot \psi_\cube$ for any axis-aligned cube $\cC$ and any measurable set $K$.
 \par
We end this section with an outline of the remainder of the paper.
In section \ref{section:prelim} we introduce definitions and conventions that will appear in subsequent sections of the paper.
In section \ref{section:cube-ub} we prove that the $\ell_0$ isoperimetric coefficient of an axis-aligned cube is of order at most $O(n^{-1/2})$.
In section \ref{section:cube-lb} we prove that the $\ell_0$ isoperimetric coefficient of an axis-aligned cube is of order at least $\Omega(n^{-1/2})$.
\section{Preliminaries}\label{section:prelim}
Throughout this paper we will use the following conventions. $A \approx B$ means that $cA < B < CA$ for some universal
constants $c$ and $C$. $B = O(A)$ means that $B < CA$ for some absolute constant $C$. $B = \Omega(A)$ means $B > cA$ for some absolute constant $c$.
$A = \Theta(B)$ means $A \approx B$. $A = o(B)$ means $ A/B \to 0 \text{ as }n \to \infty$.
 We will use $n$ to denote the ambient dimension. We will use $Q := [0,1]^n$ to denote the unit cube and $H := \{0,1\}^n$ to denote the discrete hypercube. We use $[n]$ to denote the set $\{1,\cdots, n\}$. We use $\vold{d}{\cdot}$ to denote $d$-dimensional Lebesgue measure and drop the subscript when $d = n$. In addition when we say a set is measurable we mean it is Borel measurable.
\begin{definition}\label{def:$l_0$}
We define the $\ell_0$ norm: $\R^n \rightarrow \Z^{\geq 0}$ as $\ell_0(x) = \#\{i \in [n]~:~x_i \neq 0\}$. In particular two sets $A$,$B$ are axis-disjoint if and only if $\ell_0(x-y) \geq 2$ for every $x \in A, y \in B$.
\end{definition}
\begin{definition}\label{def:cover}
Given sets $A,B$ with $A \subseteq B$ we define the $\ell_0$ boundary of $A$ with respect to $B$, which we denote as $\lobound_B{A} := \lobound A$, as the set of all points $x$ in $B \setminus A$ for which $\{x\}$ and $A$ are not axis-disjoint. Equivalently 
\[
\lobound A := \{x \in B\setminus A~:~\ell_0(x-y) \leq 1~\text{ for some }y \in A\}.
\]
\end{definition}

\begin{definition}\label{def:rectangle}
Given $x,y \in \R^n$ with $x_i \leq y_i$ for all $i \in [n]$ we write $R(x,y)$ to denote the closed rectangle with antipodal corners $x$ and $y$. Formally
\[
R(x,y) := \{z \in \R^n~:~x_i \leq z_i \leq y_i~\text{ for all }i \in [n]\}.
\]
\end{definition}
\begin{definition}\label{def:anchored}
Given $S \subseteq Q$ we say that $S$ is ``anchored'' if for any point $x \in S$ the rectangle $R(0,x)$ is a subset of $S$ (here $0$ denotes the zero vector).
\end{definition}
\begin{remark}
    In \cite{harper1999isoperimetric} the set $S$ defined in \refD{def:anchored} would be referred to as a lower set (see Appendix A from \cite{harper1999isoperimetric}), being the downward closed set in the coordinate wise partial order of vectors in $Q$.
\end{remark}
\begin{definition}\label{def:Hamm-ball}
Given $p \in [0,1]$ and $r \in \{0\} \cup [n]$ we denote by $H(p,r)$ the $p$-weighted hamming ball of radius $r$. Formally
\[
H(p,r) := \{ x \in Q~:~\#\{i : 0 \leq x_i \leq p\} \geq n-r\}.
\]
\end{definition}
\begin{definition}\label{def:grid}
Let $\{T_i\}_{i = 1}^n := \{(t_{i,j})_{j = 0}^{k_i}\}_{i = 1}^n$  denote a set of $n$ sequences where $0 = t_{i,0}, 1 = t_{i,k_i}$ for all $i \in [n]$ and $t_{i,j} < t_{i,j+1}$ for all $i \in [n], 0 \leq j < k_i$. We refer to such a set as a grid. 
\end{definition}
\begin{definition}\label{def:grid-struct}
Given a set $A \subseteq Q$ and a grid $\{T_i\}_{i = 1}^n$ we say that $A$ has a grid structure with respect to $T$ if $A$ can be written as a union of grid blocks of $\{T_i\}_{i = 1}^n$. Formally, $A$ has a grid structure with respect to $T$ if there exists a collection of $n$-tupples $\Lambda \subseteq \{(\alpha_1,\cdots,\alpha_n)~:~ 1 \leq \alpha_i < k_i\}$ such that 
\[
A = \bigcup_{\alpha \in \Lambda} [t_{1,\alpha_1},t_{1,\alpha_1 + 1}] \times \cdots \times [t_{n,\alpha_n}, t_{n,\alpha_n + 1}].
\]
Furthermore we say that $A$ has a grid structure if there exists some grid $T$ such that $A$ has a grid structure with respect to $T$.
\end{definition}
\section{Upper bound on $\psi_{\cube}$}\label{section:cube-ub}
Here we show that $\psi_{\cube} = O(n^{-1/2})$. We do this by proving a combinatorial analogue on the discrete hypercube and pass to the continuous hypercube.  We shall need a few technical lemmas. 
\begin{lemma}\label{lem:sterling}
\[
\lim_{n \to \infty}\binom{n}{\floor{n/2}}\cdot \frac{\sqrt{n}}{2^n} = \lim_{n \to \infty} \binom{n}{\ceil{n/2}}\cdot \frac{\sqrt{n}}{2^n} = \sqrt{\frac{2}{\pi}}.
\]
\end{lemma}
\begin{lemma}\label{lem:discrete-$l_0$}
There exists a choice of axis-disjoint sets $S_1,S_2 \subseteq H$ such that $S_3 := H \setminus (S_1 \cup S_2)$ satisfies $|S_3|/(\min\{|S_1|,|S_2|\}) = O(n^{-1/2})$
\end{lemma} 
\begin{proof}
Define:
\[
S_1 := \{ x \in H  ~:~ \norm{x}_1 > \ceil{n/2}\} ~~S_2 := \{ x \in H ~:~ \norm{x}_1 < \floor{n/2}\}.
\]
By symmetry $|S_1| = |S_2|$. Furthermore for every $x \in S_1, y \in S_2$ we have that $\ell_0(x-y) \geq \ell_0(x) - \ell_0(y) \geq 2$. Therefore by \refD{def:$l_0$} $S_1$ and $S_2$ are axis-disjoint. In addition $S_3 = \{ x \in H ~|~ \floor{n/2} \leq \norm{x}_1 \leq \ceil{n/2}\}$. Therefore
\[
|S_3| = \sum_{k = \floor{n/2}}^{\ceil{n/2}}\binom{n}{k}.
\]
If $n$ is even then by \refL{lem:sterling} $|S_3| = (1 + o(1))(2^n/\sqrt{n}) \cdot \sqrt{2/\pi}$. If $n$ is odd then by \refL{lem:sterling} $|S_3| = (1 + o(1))(2^{n}/\sqrt{n}) \cdot \sqrt{8/\pi}$. In either case we have that 
\[
\lim_{n \to \infty} \sqrt{n}\cdot \frac{|S_3|}{\min\{|S_1|,|S_2|\}} \leq \lim_{n \to \infty} \frac{2^n\sqrt{\frac{8}{\pi}}}{2^{n-1} - \frac{2^{n}}{\sqrt{n}}\sqrt{\frac{2}{\pi}}}  = 4\sqrt{\frac{2}{\pi}}.
\]
The result follows.
\end{proof}
\begin{corollary}\label{cor:$psi_0$-ub}
$\psi_{\cube} = O(n^{-1/2})$.
\begin{proof}
Let $S_1,S_2,S_3 \subseteq H$ be the sets guaranteed by \refL{lem:discrete-$l_0$}. 
\\Define $I_0 := [0,1/2), I_1 := [1/2,1]$ and the map $f : H \rightarrow \mathcal{P}(Q)$ according to 
\[
f(x) := I_{x_1} \times \cdots \times I_{x_n}.
\]
Note that $f$ maps every point in $H$ to some cube of side length 1/2 contained in $Q$. We now record some properties of $f$:
\begin{itemize}
\item $\vol{f(x)} = 2^{-n}~\text{ for all } x \in H$: This follows that from the fact that $f(x)$ is the cartesian product of $n$ intervals, each of length $1/2$.
\item $x \neq y \implies f(x) \cap f(y) = \emptyset$: This follows from the fact that $x \neq y$ means $x_i \neq y_i$ for index $i \in [n]$. Since $I_0$ and $I_1$ are disjoint we conclude that $I_{x_i} \cap I_{y_i} = 0$ and the cubes which $x$ and $y$ map to are disjoint.
\item $\cup_{x \in H}f(x) = Q$: Since $I_0 \cup I_1 = [0,1]$ every point in $f(x)$ has coordinates values in $[0,1]$ and is therefore in $Q$. On the other hand every point $p \in Q$ is contained in the cube $I_{\min\{1,\floor{2p_1}\}} \times \cdots \times I_{\min\{1,\floor{2p_n}\}}$. Since $(\min\{1,\floor{2p_1}\},\cdots,\min\{1,\floor{2p_n}\}) \in H$ the result follows.
\item If $x,y \in H$ are axis-disjoint then $f(x),f(y)$ are axis disjoint: This follows from the fact that $\ell_0(x-y)_1 > 1$ means $x$ and $y$ disagree on at least 2 coordinates (say $i_1,i_2$). Therefore $I_{x_{i_1}} \cap I_{y_{i_1}} = I_{x_{i_2}} \cap I_{y_{i_2}} = \emptyset$ and every point in $f(x)$ disagrees with every point in $f(y)$ on at least 2 coordinates. 
\end{itemize}
As a consequence $f(S_1),f(S_2)$ are axis-disjoint, $f(S_3) = Q\setminus(f(S_1) \cup f(S_2))$ and by \refL{lem:discrete-$l_0$}
\[
\frac{\vol{f(S_3)}}{\min\{\vol{f(S_1)},\vol{f(S_2)}\}} = \frac{|S_3|}{\min\{|S_1|,|S_2|\}} = O(n^{-1/2}).
\] 
\end{proof}
\end{corollary}
\section{Lower bound on $\psi_{\cube}$}\label{section:cube-lb}
\subsection{Proof Strategy}
As defined in \cite{laddha2020convergence}, $\psi_{\cube}$ is the minimum value of $\vol{Q \setminus (S_1 \cup S_2)}/\min\{\vol{S_1},\vol{S_2}\}$ across all axis-disjoint subsets $S_1,S_2 \subseteq Q$. Note that this definition assumes $\cC = Q$ since the isoperimetric coefficient of an axis aligned cube is invariant under translation in any direction and dilation. Note that if $S_1,S_2$ are axis-disjoint then $\lobound(S_1) \subseteq Q\setminus(S_1 \cup S_2)$. Our proof strategy is motivated by the following observation 
\begin{proposition}\label{prop:coeff-l0}
    \begin{equation}\label{eq:iso-coeff-equiv}
        \inf_{\underset{S_1,S_2 - \text{axis-disjoint}}{S_1,S_2 \subset K}} \frac{\vol{K\setminus(S_1 \cup S_2)}}{\min\{\vol{S_1},\vol{S_2}\}} =  \inf_{\underset{S_1,S_2 - \text{axis-disjoint}}{S_1,S_2 \subset K}} \frac{\vol{\lobound(S_1)}}{\min\{\vol{S_1},\vol{S_2}\}}.
    \end{equation}
\end{proposition}
\begin{proof}
    The left hand side is at most the right hand side since $\lobound(S_1) \subset K\setminus(S_1\cup S_2)$ for every pair of axis-disjoint sets $S_1,S_2$. On the other hand, given $S_1,S_2$ the set $S_2' := S_2 \cup (K\setminus(S_1 \cup S_2)) \setminus \lobound(S_1)$ is axis-disjoint from $S_1$, has volume at least as large as $S_2$, and $K\setminus(S_1 \cup S_2') = \lobound(S_1)$. Since both sides are infimums over pairs of axis-disjoint sets it follows that the right hand side is at most the left hand side.
\end{proof}

\refP{prop:coeff-l0} says that lower bounding $\psi_\cC$ reduces to showing that for any pair of axis-disjoint sets that the $\ell_0$ boundary of one of the two sets is large. Thus it is in our interest to know what types of axis-disjoint sets minimize the left-hand side of the desired inequality.
To that end we state the following isoperimetric theorem due to Harper \cite[Theorem 1 + Theorem 2]{harper1999isoperimetric}.
\begin{theorem}\label{thm:Harper-grid}
Let $\lambda \in (0,1)$ and let $\mathcal{A}$ be the collection of all subsets of $Q$ that are anchored and have a grid structure. Then there exists a choice of $0 < p < 1$ and $0 \leq r \leq n$ for which
\begin{equation}\label{eq:Harper}
\min_{\underset{\vol{A} = \lambda}{A \in \mathcal{A}}}\vol{\lobound A} = \vol{H(p,r)},
\end{equation}
where $H(p,r)$ is a $p$-weighted hamming ball (see \refD{def:Hamm-ball}) satisfying $\vol{H(p,r)} = \lambda$. 
\end{theorem}
\begin{remark}\label{rem:Harper-grid}
We note that Theorem 1 + Theorem 2 of \cite{harper1999isoperimetric} only directly implies \refT{thm:Harper-grid} when we further restrict $\cA$ to subsets with grids of the form $\{\{0,t_i,1\}\}_{i = 1}^n$. However, Harper argued that the minimizers in $\cA$ are exactly of this more restricted form (see pages 288 - 289 of \cite{harper1999isoperimetric}, in particular bullet point 2 of page 289). For such sets Theorem 1 of \cite{harper1999isoperimetric} then tells us that the left hand side of \eqref{eq:Harper} is obtained by a set of the form $\{x \in Q~:~\#\{i:0 \le x_i \le t_i\} \ge n - r\}$ for some $r \in [n]$. Theorem 2 of \cite{harper1999isoperimetric} then tells us that this set should satisfy $t_i = t_j$ for all $i,j \in [n]$. In particular the set is equal to $H(p,r)$ for some $p \in (0,1), r \in [n]$. 
\end{remark}
In addition we note that $\lobound H(p,r) = H(p,r+1)\setminus H(p,r)$.
Therefore if $S_1$ was anchored and had a grid structure we could apply \refT{thm:Harper-grid} with $A := S_1$ and $\lambda = \vol{S_1}$ to reduce the problem to showing that 
\[
\inf_{p,r} \frac{\vol{H(p,r+1)\setminus H(p,r)}}{\min\{\vol{H(p,r),\vol{Q\setminus H(p,r+1)}}\}} \geq \Omega(n^{-1/2}).
\]
Of course $S_1$ need not be anchored or have a grid structure so we need a procedure to transform $S_1$ into a set of the desired type while maintaining axis-disjointness between $S_1$ and $S_2$ and preserving the volumes of $S_1$ and $S_2$. In particular we need the following lemmas:
\begin{lemma}\label{lem:step 1}
Let $S_1,S_2 \subseteq Q$ be axis-disjoint. Then there exists axis-disjoint sets $S_1',S_2' \subseteq Q$ such that $\vol{S_1'} = \vol{S_1}, \vol{S_2'} = \vol{S_2}$, 
and $S_1'$ is compact, anchored (see \refD{def:anchored}), and has a boundary of measure 0.
\end{lemma}
\begin{lemma}\label{lem:step 2}
Let $S_1,S_2 \subseteq Q$ be axis-disjoint with $S_1$ being anchored, compact, and having a boundary of measure 0. Then for every $\eps > 0$ there exists axis-disjoint subsets $S_1',S_2' \subseteq Q$ such that $S_1'$ is anchored, has a grid structure and satisfies 
\[
\frac{\vol{\lobound(S_1)}}{\min\{\vol{S_1},\vol{S_2}\}} \geq \frac{\vol{\lobound(S_1')}}{\min\{\vol{S_1'},\vol{S_2'}\}} - \eps.
\]
\end{lemma}
\begin{lemma}\label{lem:step 3}
There exists an absolute constant $c > 0$ such that
\[
\inf_{p,r} \frac{\vol{H(p,r+1)\setminus H(p,r)}}{\min\{\vol{H(p,r),\vol{Q\setminus H(p,r+1)}}\}} \geq cn^{-1/2}.
\]
\end{lemma}
We quickly prove \refT{thm:cube} using these lemmas:
\begin{proof}[Proof of \refT{thm:cube}]
By \refC{cor:$psi_0$-ub} it suffices to show that $\psi_{\cube} = \Omega(n^{-1/2})$. Let $S_1,S_2 \subseteq Q$ be arbitrary axis-disjoint subsets. By applying \refP{prop:coeff-l0}, \refL{lem:step 1}, \refL{lem:step 2} with $\eps \to 0^+$, \refT{thm:Harper-grid} and \refL{lem:step 3} we conclude that
\[
\frac{\vol{Q\setminus(S_1\cup S_2)}}{\min\{\vol{S_1},\vol{S_2}\}} \geq \inf_{p,r} \frac{\vol{H(p,r+1)\setminus H(p,r)}}{\min\{\vol{H(p,r),\vol{Q\setminus H(p,r+1)}}\}} \geq \Omega(n^{-1/2}).
\]
Since $S_1,S_2$ were arbitrary the result follows.
\end{proof}
The proof of \refL{lem:step 2} is straightforward while the proofs of \refL{lem:step 1} and \refL{lem:step 3} are more technically challenging.
\subsection{Proof of \refL{lem:step 1}: Shaking} 
Given $S_1,S_2$ we will construct $S_1',S_2'$ satisfying the guarantee of \refL{lem:step 1} using the so called shaking operation.
A known technique in the study of vertex-isoperimetry (in the form of compression) and convex geometry, shaking is the process of taking a set $X \subseteq Y$, a direction $v \in S^{n-1}$ and ``shaking" every fiber of $X$ along the $v$ direction into an initial segment of $Y$ along the $v$ direction. With this in mind our plan is to shake $S_1$ towards the zero vector and $S_2$ towards the all-ones vector, while maintaining axis-disjointness and volume. We will also show that at the end of the process $S_1'$ is compact, anchored, and has a boundary of measure 0.
\par
Here we define what it means to ``shake" axis disjoint sets along a coordinate axis-direction. 
\begin{definition}\label{def:shake_i}
Let $i \in [n]$. We define the shaking transforms $\shake^+_i,\shake^-_i$ as follows: Let $X \subset Q$ be measurable. Then
\begin{align}\label{eq:Shaking}
\shake^+_i(X) &:= \bigcup_{x \in \Proj_{e_i^\top} X} \{x\} \times \left[ 0, \vold{1}{\{ y \in X~|~\Proj_{e_i^\top}y = x \}} \right], \\
\shake^-_i(X) &:= \bigcup_{x \in \Proj_{e_i^\top} X} \{x\} \times \left( 1 - \vold{1}{\{ y \in X~|~\Proj_{e_i^\top}y = x \}},1 \right),
\end{align}
where $e_i$ denotes the $i$th standard basis vector. Furthermore we define 
\begin{align}\label{eq:Shaking-full}
\shake^+ := \shake^+_n \circ \shake^+_{n-1} \circ \cdots \circ \shake^+_1, \\
 \shake^- := \shake^-_n \circ \shake^-_{n-1} \circ \cdots \circ \shake^-_1.
\end{align}
\end{definition}
Here we record some of the properties of the transformed sets:
\begin{lemma}\label{lem:transform-prop}
Let $A,B \subseteq Q$ be measurable axis-disjoint sets. Then $(\shake^+_i(A),\shake^-_i(B)) := (\hat{A},\hat{B})$ satisfies the following:
\begin{enumerate}
\item $\hat{A},\hat{B}$ are measurable.
\item $\vol{A}= \vol{\hat{A}},\vol{B}= \vol{\hat{B}}$.
\item $\hat{A},\hat{B}$ are axis-disjoint.
\item $\vol{\lobound\hat{A}} \leq \vol{\lobound A}, \vol{\lobound\hat{B}} \leq \vol{\lobound B}$.
\end{enumerate}
\end{lemma}
The proof of this lemma is included in the appendix. Using the transformation from \refD{def:shake_i} our desired sets are $S_1' := \shake^+(\shake^+(S_1))$ and $S_2' := \shake^-(\shake^-(S_2))$.
\refL{lem:transform-prop} guarantees that $S_1',S_2'$ have the desired volumes and are axis-disjoint. What is left to show is that $S_1'$ is compact and has a boundary of measure 0. We first prove compactness.
\begin{lemma}\label{lem:closed-rec}
Let $A\subseteq Q$ be measurable. Let $\hat{A} = \shake^+(A)$ and let $x \in \hat{A}$. Then $R(0,x)$ is contained in $\hat{A}$.
\begin{proof}
We prove this by induction, showing that for $\hat{A}_k := \shake^+_k(\shake^+_{k-1}(\cdots (\shake^+_1(A)\cdots))$, if $x$ is in $\hat{A}_k$ the rectangle  $[0,x_1] \times \cdots \times [0,x_k] \times \{x_{k+1}, \cdots,x_n\}$ is contained in $\hat{A}_k$. In the base case, $k = 1$, the result is immediate from the definition of shaking as every fiber of $A$ in the direction $e_1$ is transformed into a closed initial segment of $Q$. Assuming the inductive hypothesis holds up to $k$, for the case of $k+1$ we argue as follows: Let $x$ be in $\hat{A}_{k+1}$. Take $Y := \{y \in \hat{A}_{k}~|~ y_i = x_i ~\forall i \neq k+1\}$. Note that $x \in \phi_{k+1}^+(Y)$ and that $x_{k+1} \leq \vold{1}{Y}$. By the induction hypothesis the rectangle $R_{y,k} := [0,x_1] \times \cdots \times [0,x_k] \times \{y_{k+1}\} \times \{x_{k+2} \cdots,x_n\}$ is contained in $\hat{A}_k$ for all $y \in Y$. Therefore every fiber $f$ in the direction of $e_{k+1}$, whose first $k$ coordinates are contained in $[0,x_1] \times \cdots \times [0,x_k]$ and whose last $n-k-1$ coordinates agree with $x$, intersects $\hat{A}_k$ on a set of 1 dimensional measure equal to $\vold{1}{Y}$. Therefore every such fiber intersects $\hat{A}_{k+1}$ on an initial segment of length $\vold{1}{Y}$. Therefore $\{z_1, \cdots, z_k\} \times [0, \vold{1}{Y}] \times \{x_{k+2},\cdots,x_n\} \subseteq \hat{A}_{k+1}$ for all $z$ satisfying $0 \leq z_i \leq x_i$ for all $1 \leq i \leq k$. Therefore  $[0,x_1] \times \cdots \times [0,x_{k+1}] \times \{x_{k+2}, \cdots,x_n\}$ is contained in $\hat{A}_{k+1}$.

\end{proof}
\end{lemma}
\begin{corollary}\label{cor:comp}
Let $A \subseteq Q$ be measurable. Let $\hat{A} = \shake^+(A)$ and $\hat{C} = \shake^+(\shake^+(A))$. Then the closure of $\hat{A}$ is $\hat{C}$.
\begin{proof}
If $p \in \hat{A}$ then by \refL{lem:closed-rec} the rectangle $R(0,p)$ is contained in $\hat{A}$.  Suppose now that $x$ is a limit point of $\hat{A}$. Then there exists a sequence of points $(p_i)_{i \geq 1}$ in $\hat{A}$ for which $\norm{p_i - x}_\infty$ is decreasing and tends to 0. Since $\norm{p_i - x}_\infty$ is decreasing we conclude that $[0,x_1) \times \cdots \times [0,x_n) \subseteq \cup_{i \geq 1} R(0,p_i) \subseteq \hat{A}$. From the definition of the shaking operation we may further deduce that $[0,x_1] \times \cdots \times [0,x_k] \times [0,x_{k+1}) \times \cdots \times [0,x_n) \subseteq \shake_k(\shake_{k-1}( \cdots \shake_1(\hat{A})\cdots)$ for all $k \geq 1$. Therefore $x \in \hat{C}$. Suppose now that $\hat{C}$ contains a point $y$ outside the closure of $\hat{A}$. Then its distance to the closure of $\hat{A}$ is positive. In particular $R(0,y)\setminus \hat{A}$ is a set of non-zero measure. But by \refL{lem:closed-rec} $R(0,y) \subset \hat{C}$ which means $\vol{\hat{C}} > \vol{\hat{A}}$, contradicting \refL{lem:transform-prop}. 
\end{proof}
\end{corollary}
By \refC{cor:comp} $S_1'$ is compact. We now prove that the boundary of $S_1'$ has 0 volume: \begin{corollary}\label{cor:meas0}
Let $A \subseteq Q$ be measurable. Then $\hat{C} := \shake^+(\shake^+(A))$ has boundary of measure 0.
\begin{proof}
We first recall the well known Steinhaus theorem \cite{stromberg1972elementary} which says that the difference set of a subset of $\R^n$ of positive measure contains a  ball, centered at the origin, of positive radius. We will show that $\partial\hat{C} - \partial\hat{C}$ does not contain such a ball. In particular we will show $\partial\hat{C} - \partial\hat{C}$ does not contain a point in the positive orthant. To this end assume for the sake of contradiction that it does. Then there exists $a,b \in \partial\hat{C}$ for which $b-a$ is in the positive orthant. In particular every entry of $a$ is strictly smaller than every entry of $b$. By \refL{lem:closed-rec} $\hat{C}$ contains the rectangle $R(0,b)$, meaning $a$ is an interior point of $R(0,b)$ and therefore $\hat{C}$. But $a \in \partial\hat{C}$  and we arrive at a contradiction.
\end{proof}
\end{corollary}
We have thus shown that $S_1',S_2'$ satisfy the guarantee of \refL{lem:step 1}.
\subsection{Proof of \refL{lem:step 2}}
\begin{proof}[Proof of \refL{lem:step 2}]
Let $S_1,S_2$ be the sets guaranteed by \refL{lem:step 1}.
 Since $S_1$ is compact, anchored, and has a boundary of measure 0, for any $\eps > 0$ there exists a choice of axis-aligned grid cubes $\mathcal{D}_\eps$ of sufficiently small side length such that $\mathcal{D}_\eps \subseteq S_1$ and $\vol{S_1\setminus \mathcal{D}_\eps} < \eps$ and $\mathcal{D}_\eps$ is anchored and has a grid-structure. A brief sketch of this is as follows: Every point in the interior of $S_1$ is some non-zero distance away from the boundary of $S_1$. Therefore there exists $s > 0$ such that the measure of the set of points in the interior of $S_1$ that are distance more that $s$ from the boundary is at most $\eps$.  
 Since the diameter of an axis-aligned cube of side-length $s$ is $s\sqrt{n}$ and the boundary of $S_1$ has measure 0, we conclude that, except for a set of measure at most $\eps$, all points in $S_1$ are contained in axis-aligned grid cubes of side length $s/\sqrt{n}$ that also lie in $S_1$. Finally, since $S_1$ is anchored, this set of axis-aligned grid cubes can be extended to a set of axis-aligned grid cubes that is also anchored and contained in $S_1$.
 \par
 Note now that
\[
\frac{\vol{\lobound(\cD_\eps)}}{\min\{\vol{\cD_\eps},\vol{S_2}\}} \leq \frac{\vol{\lobound(S_1)} + \eps}{\min\{\vol{S_1} - \eps,\vol{S_2}\}} \leq \frac{\vol{\lobound(S_1)}}{\min\{\vol{S_1},\vol{S_2}\}} + \eps',
\]
where $\eps' \to 0$ as $\eps \to 0$. \end{proof}
\subsection{Proof of \refL{lem:step 3}: An Isoperimetric Inequality for Binomials}
Because $Q$ has volume 1, by the definition of $H(p,r)$ (\refD{def:Hamm-ball}), one may view the $\ell_0$-isoperimetry lowerbound of $p$-weighted hamming balls as an isoperimetric inequality about binomial random variables. This is because for any $k \in [n]$, $\vol{H(p,k)}$ is the probability that an $(n,1-p)$ binomial random variable is at most $k$, and  $\vol{H(p,k)\setminus H(p,k-1)}$ is the probability that such a random variable is equal to $k$. With this connection in mind we state the following inequality for binomial random variables: 
\begin{lemma}\label{lem:binom-bound}
There exists an absolute constant $c > 0$ such that the following is true. Let $X$ be distributed as an $(n,1-p)$ binomial random variable. Then for all $k \in \{0,1,\cdots,n\}$ the following is true:
\begin{equation}\label{eq:binom}
\bP[X = k] \geq \frac{c}{\sqrt{np(1-p)}}\min\{ \bP[X < k],\bP[X > k] \}.
\end{equation}
\end{lemma}
Although \refL{lem:binom-bound} is known we provide a proof of it in the appendix for completeness. Note that when proving \refL{lem:binom-bound} we may assume that $p \leq 1/2$. This is because for $X$ an $(n,1-p)$ binomial random variable and $\tilde{X}$ an $(n,p)$ binomial random variable, one has 
\[
\min\{\bP[X > k],\bP[X < k]\} = 
\min\{\bP[\tilde{X} < n-k],\bP[\tilde{X} > n-k]\},
\]
for every $k \in \{0,1,\cdots,n\}$, and the inequality guaranteed by \refL{lem:binom-bound} holds for all $k \in \{0,1,\cdots,n\}$. 
\subsection{A generalization of \refT{thm:cube}}
The following observation was made by Vladimir Koltchinskii after being presented \refT{thm:cube}. 
\begin{corollary}\label{cor:cube}
    Let $\Omega = \Omega_1 \times \cdots \times \Omega_n$, where $(\Omega_j, {\mathcal A}_j, P_j)$ are probability spaces. Suppose there exist measurable bijections $\phi_j :[0,1]\mapsto \Omega_j$ such that 
$P_j=\lambda\circ \phi_j^{-1},$ where $\lambda$ is the Lebesgue measure. Let $P=P_1\times \dots \times P_n$. Define 
\[
\psi_\Omega := \inf_{\underset{S_1,S_2 -\text{axis-disjoint}}{S_1,S_2 \subset \Omega}} \frac{P(\Omega\setminus(S_1 \cup S_2))}{\min\{P(S_1),P(S_2)\}}
\]
Then $\psi_\Omega = \Theta(n^{-1/2})$.
\end{corollary}
The corollary is obtained by observing that $S_1,S_2 \subset \Omega$ are axis-disjoint if and only if $\phi^{-1}(S_1),\phi^{-1}(S_2)$ are axis-disjoint in $Q$ and that $P(S) = \lambda^n(\phi^{-1}(S))$. The result is then an immediate consequence of \refT{thm:cube}. A natural setting for which the corollary holds is when $\Omega_j$ is a Polish space, $\cA_j$ is a Borel $\sigma$-algebra and $P_j$ is a non-atomic probability measure. 
\section{Upper bounds for General Bodies}\label{section:bodies-ub}
Given the result of the previous section is it is natural to ask if the cube is ``special" in terms of $\ell_0$-isoperimetry. Remarkably, among \textit{all} measurable sets the $\ell_0$-isoperimetric coefficient of a cube is, up to an absolute constant, largest.  
\begin{theorem}\label{thm:uni-ub}
If $K \subseteq \R^n$ is measurable then $\psi_K = O(n^{-1/2})$. 
\end{theorem}
We remark here that in the case where $n$ is not too large we will use an ad-hoc argument. Therefore, up until the proof of \refT{thm:uni-ub}, we will assume that $n \geq n_0$ for $n_0$ a sufficiently large constant.
Our proof of \refT{thm:uni-ub} relies on extending our approach for upperbounding the $\ell_0$-isoperimetric coefficient of the cube to general bodies. As done in the cube case we will assume that $\vol{K} = 1$ (scaling and translation does not effect the $\ell_0$ isoperimetric coefficient of a set). 
Extending the upper bound approach is relatively straightforward:
\begin{enumerate}
\item Given a measurable set $K \subseteq \R^n$ take $p \in \R^n$ to be any point satisfying 
\[
\vol{K \cap \{x_i \leq p_i\}} = \vol{K \cap \{x_i > p_i\}}, 
\]
 for all $i \in [n]$. Note that such a $p$ always exists. To see this recall that Lebesgue measure is not supported on any $n-1$ dimensional hyperplane. Therefore the map $t \mapsto \vol{K \cap \{x_i \leq t\}}$ is continuous for every $i \in [n]$ and its image contains $(0,1)$. By the intermediate value theorem there exists some $p_i$ for which $\vol{K \cap \{x_i \leq p_i\}} = 1/2 = \vol{K \cap \{x_i > p_i\}}$. 
\item Write $K = \cup_{s \in H} K_s$ where for every $x \in K$ we have
\[
x \in K_s \iff 
\begin{cases}
x_i \leq p_i & \text{ if } s_i = 0\\
x_i > p_i & \text{ if } s_i = 1
\end{cases}.
\]
\item Take $S_1,S_2$ to be axis-disjoint subsets of $H$ and $S_3 := H \setminus (S_1 \cup S_2)$.
\item Observe that 
\[
\psi_K \leq \frac{\vol{\cup_{s \in S_3}K_s}}{\min\left\{\vol{\cup_{s \in S_1} K_s},\vol{\cup_{s \in S_2} K_s}\right\}}.
\]
\end{enumerate}
In the case of the hypercube we took $p := (0.5,\cdots,0.5)$, $S_1 := \cup_{\# \{i:s_i = 1\} < \floor{n/2}}K_s$ and $S_2 := \cup_{\#\{i:s_i = 1\} > \ceil{n/2}}K_s$.  This particular construction can be immediately generalized to bodies for which $\vol{K_s}$ is the same for all $s$, giving the same upperbound of $O(n^{-1/2})$ for the isoperimetric coefficient. A general class of such bodies are unconditional bodies, which are measurable compact bodies that are closed under reflections over coordinates hyperplanes.
\begin{corollary}\label{cor:uncond-body-ub}
Let $K$ be an unconditional body. Then $\psi_K = O(n^{-1/2})$.
\end{corollary}
In general however it is not clear how to pick $S_1, S_2$, as it can be difficult to determine axis-disjoint sets, each with a constant fraction of the mass, whenever the distribution of the mass amongst the orthants is non-uniform. 
\subsection{Constructing $S_1,S_2$ with a random splitting plane.}
\sloppy 
To implement the upper bound approach for more general bodies, we will show that one can construct $S_1,S_2,S_3$ randomly. In particular, we will show that a ``random splitting plane" has non-zero probability of separating $H$ into 3 pieces $S_1,S_2,S_3$ where $S_1,S_2$ are axis-disjoint, $\vol{\cup_{s \in S_3}K_s} = O(n^{-1/2})$ and $\min\{\vol{\cup_{s \in S_1}K_s},\vol{\cup_{s \in S_2}K_s}\} = \Omega(1)$.
We start with a definition.
\begin{definition}\label{def:balanced}
A function $w : H \to [0,1]$ is ``balanced'' if it satisfies the following properties:
\begin{enumerate}
\item $\sum_{s \in H} w(s) = 1$.
\item $\sum_{\underset{s_i = 0}{s \in H}}w(s) = \sum_{\underset{s_i = 1}{s \in H}}w(s)$ for all $i \in [n]$.
\end{enumerate}
\end{definition}
Note that every body $K$ of volume 1 gives rise to a ``balanced'' function $w$, where $w(s) := \vol{K_s}$.
Given $x,y \in H$ we write $\ip{x}{y} := \#\{i~s.t.~ x_i \neq y_i\}$. Given $s \in H$ we write $\bar{s} = (1,1,\cdots,1) - s$.
We now give our construction of $S_1$ and $S_2$.
\begin{enumerate}
\item Pick $z$ uniformly from $H$.
\item Take
\[
S_1 := \{s \in H~:~\ip{s}{z} < \floor{n/2}\}, S_2 := \{s \in H~:~\ip{s}{z} > \ceil{n/2}\}, S_3 := \{s \in H~:~\ip{s}{z} \in [\floor{n/2},\ceil{n/2}]\}.
\]
\end{enumerate} 
Note that $S_1,S_2$ are axis-disjoint since for $s_1 \in S_1, s_2 \in S_2, \ip{s_1}{s_2} \geq \ip{s_1}{z}-\ip{s_2}{z} \geq 2$, meaning $s_1,s_2$ disagree on at least 2 coordinates.
Since $s \in S_1 \iff \ip{s}{z} < \floor{n/2}$ we have that
\begin{align*}
\bE[\vol{\cup_{s \in S_1} K_s}] = \sum_{s \in H}w(s)\cdot \bP(\ip{s}{z} < \floor{n/2}) 
&= \bP(\ip{(0,0,\cdots,0)}{z} < \floor{n/2}) \cdot \sum_{s \in H}w(s),\\
&= \bP(\ip{(0,0,\cdots,0)}{z} < \ceil{n/2}), \\
&= \frac{\sum_{k = 0}^{\floor{n/2} - 1}\binom{n}{k}}{2^n}.
\end{align*}
Similar computations show that $
\bE[\vol{\cup_{s \in S_2} K_s}] =  \frac{\sum_{k = \ceil{n/2} + 1}^{n}\binom{n}{k}}{2^n}$ and $\bE[\vol{\cup_{s \in S_3} K_s}] =   \frac{\sum_{k = \floor{n/2}}^{\ceil{n/2}}\binom{n}{k}}{2^n}$.
By \refL{lem:sterling}  $\bE[\vol{\cup_{s \in S_3} K_s}]  = O(n^{-1/2})$. By Markov's inequality we conclude, with any constant probability less than 1, that $\vol{\cup_{s \in S_3} K_s}  = O(n^{-1/2})$. We are left with dealing with $S_1$ and $S_2$. When $K$ is symmetric $\min\{\vol{\cup_{s \in S_1} K_s},\vol{\cup_{s \in S_2} K_s}\}$ is essentially 1/2.
\begin{lemma}\label{lem:symm}
Let $K \subseteq \R^n$ be symmetric. Then for every $z \in H$ we have that 
then 
\[
\sum_{\ip{s}{z} < \floor{n/2}}w(s) = 
\sum_{\ip{s}{z} > \ceil{n/2}}w(s).
\]
\end{lemma}
\begin{proof}
Since $K$ is symmetric $w(s) = w(\bar{s})$ for every $s \in H$. Furthermore $\ip{s}{z} < \floor{n/2} \implies \ip{\bar{s}}{z} >  n - \floor{n/2} = \ceil{n/2}$. Since $s \mapsto \bar{s}$ is a bijection on $H$ the result follows.
\end{proof}
\refL{lem:symm} implies that $\min\{\vol{\cup_{s \in S_1} K_s},\vol{\cup_{s \in S_2} K_s}\} = 0.5 - \vol{\cup_{s \in S_3} K_s}$.
Thus \refL{lem:symm} and Markov's inequality imply that with non-zero probability there exists axis-disjoint subsets $S_1,S_2$ and complement $S_3$ for which $\min\{\vol{\cup_{s \in S_1} K_s},\vol{\cup_{s \in S_2} K_s}\} = 0.5 - O(n^{-1/2}) = \Omega(1)$ and $\vol{\cup_{s \in S_3} K_s} = O(n^{-1/2})$. Since the probability is non-zero, by the Probabilistic Method there exists a realization of $S_1,S_2,S_3$ with these guarantees. In particular $\psi_K = O(n^{-1/2})$.
\begin{corollary}\label{cor:symm}
Let $K$ be a symmetric convex body. Then $\psi_K = O(n^{-1/2})$.
\end{corollary}
\subsection{Bounding the variance of $\vol{S_1}$}
Note that \refC{cor:symm} required the use of \refL{lem:symm}, which is not applicable when $K$ is non-symmetric.
To apply the splitting plane method for more general bodies one must show two-sided concentration of $\vol{\cup_{s \in S_1} K_s}$ (i.e. bound the variance of $\vol{\cup_{s \in S_1} K_s}$). This will imply a sufficiently good lowerbound on $\vol{\cup_{s \in S_2} K_s}$ since $\vol{\cup_{s \in S_3} K_s} = O(n^{-1/2})$ with sufficiently small constant failure probability. Let $z$ denote the random vector distributed uniformly on $H$.
Letting $I_{s}$ denote the indicator of the event $\{\ip{s}{z} < \floor{n/2}\}$ we may write
 
\begin{align*}
\bE[\vol{\cup_{s \in S_1} K_s}^2] = \sum_{s, t \in H} \vol{K_s}\vol{K_t}\bE[I_sI_t] =  \sum_{s \in H}\sum_{k = 0}^n\sum_{\underset{\ip{s}{t} = k}{t \in H}} \vol{K_s}\vol{K_t} \bE[I_sI_t].
\end{align*}
We will exploit this representation in a few ways. Firstly, if $k = \ip{s}{t}$ then (up to permutation of indices) we may write $s = (u,v), t = (u,\bar{v})$, where the length of $u$ is $n-k$ and the length of $v$ is $k$. Since $z$ is distributed uniformly over $H$, $\bE[I_sI_t]$ does not change when the same set of coordinates in $s$ and $t$ are both flipped. Thus we may assume that $u = \textbf{0}_{n-k}, v = \textbf{0}_k$, and $\bar{v} = \textbf{1}_k$. Under these assumptions, $I_sI_t$ is equal to 1 exactly when $\norm{z}_1 < \floor{n/2}$ and $\norm{z_1}_1 + k - \norm{z_2}_1 < \floor{n/2}$, where $z = (z_1,z_2)$ and $z_2$ has length $k$.  
\begin{lemma}\label{lem:decreasing}
Let $s_1,s_2,t_1,t_2 \in H$. If $\ip{s_1}{t_1} \leq \ip{s_2}{t_2}$ then 
\[
\bE[I_{s_1}I_{t_1}] \geq \bE[I_{s_2}I_{t_2}].
\]
\end{lemma}
\begin{proof}
Let $k_1 := \ip{s_1}{t_1},k_2 := \ip{s_2}{t_2}$. When $k_1 = k_2$ the result is immediate by the preceding observation. For $k_1 < k_2$ by transitivity it suffices to consider the case of $k_2 = k_1 + 1$. By the preceding observation we may assume that $s_1 = s_2 = \textbf{0}_n, t_1 = \textbf{0}_{n-k_1}\textbf{1}_{k_1}$ and $t_2 = \textbf{0}_{n-k_1-1}\textbf{1}_{k_1+1}$. Now write $z = (z_1,z_2,z_3)$ where the length of $z_1$ is $n - k_1 -1$, the length of $z_2$ is 1, and the length of $z_3$ is $k_1$. From the definition of the indicator it follows that
\begin{align*}
I_{s_1}I_{t_1} = 1 &\iff \norm{z}_1 < \floor{n/2}, \norm{z_1}_1 + \norm{z_2}_1 + k_1 - \norm{z_3} < \floor{n/2}.
\\
I_{s_2}I_{t_2} = 1 &\iff \norm{z}_1 < \floor{n/2}, \norm{z_1}_1 +1 - \norm{z_2}_1 + k_1 - \norm{z_3} < \floor{n/2}.
\end{align*}

We will now show that the number of choices of $z$ for which $I_{s_1}I_{t_1} = 1$ is at least as large as the number of choices of $z$ for which $I_{s_2}I_{t_2} = 1$.
To that end suppose $z$ satisfies $I_{s_2}I_{t_2} = 1$. We have two cases:
\begin{enumerate}
    \item $z_2 = 0$: In this case 
    \[
    \norm{z_1}_1 + \norm{z_2}_1 + k_1 - \norm{z_3}_1 < \norm{z_1}_1 + 1 - \norm{z_2}_1 + k_1 - \norm{z_3} < \floor{n/2}.
    \]
    Thus $I_{s_1}I_{t_1} = 1$.
    \item $z_2 = 1$: If $I_{s_1}I_{t_1} = 1$ then we are done. If $I_{s_1}I_{t_1} = 0$ then $\norm{z_1}_1 + 1 - \norm{z_2}_1 + k_1 - \norm{z_3}_1 < \floor{n/2}$ while $\norm{z_1}_1 + \norm{z_2}_1 + k_1 - \norm{z_3}_1 \ge \floor{n/2}$.
    Consider now the map $z \mapsto \tilde{z}$ that flips the $(n-k_1)$th entry of $z$ and fixes all other. This is a bijection on $H$. Furthermore $\norm{\tilde{z}}_1 < \norm{z}_1 < \floor{n/2}$. On the other hand $\norm{\tilde{z}_1}_1 + 1 - \norm{\tilde{z}_2}_1 + k_1 - \norm{\tilde{z}_3}_1 \ge \floor{n/2}$ while $\norm{\tilde{z}_1}_1 + \norm{\tilde{z}_2}_1 + k_1 - \norm{\tilde{z}_3}_1 < \floor{n/2}$. Therefore $\tilde{z}$ is a choice for which $I_{s_1}I_{t_1} = 1$ while $I_{s_2}I_{t_2} = 0$.
\end{enumerate}
From these two cases we have shown that the number of choices of $z$ for which $I_{s_1}I_{t_1} = 1$ is at least as large as the number of choices of $z$ for which $I_{s_2}I_{t_2} = 1$. Since $z$ is distributed uniformly over $H$ we conclude that 
\[
\bE[I_{s_1}I_{t_1}] \geq \bE[I_{s_2}I_{t_2}].
\]
\end{proof}
The previous lemma says that the probability of two vectors both ending up in $S_1$ is a decreasing function of the $\ell_1$ distance between them. We now determine a value $k$ for which $\bE[I_sI_t]$ is separated from 1/2, by an absolute constant, whenever $\norm{s-t}_1 \geq k$. 
\begin{lemma}\label{lem:fluctuation}
Let $z = (z_1,z_2)$ be distributed uniformly over $H$, where $z_2$ is of length $\floor{n/3}$. Define the following events: 
\begin{align*}
\cE &:= \{ \norm{z}_1 < \floor{n/2} \}, \\
\cF &:= \{ \norm{z_1}_1 + \floor{n/3}-\norm{z_2}_1 < \floor{n/2}\}.
\end{align*}
Then there exists an absolute constant $c > 0$ for which 
\[
\bP(\cE \cap \cF) \leq \frac{1}{2} - c.
\]\end{lemma}
The proof of this lemma is included in the appendix. 
As a corollary $\bE[I_sI_t] \leq 1/2 - c$ for $s = \textbf{0}_n, t = \textbf{0}_{n - \floor{n/3}}\textbf{1}_{\floor{n/3}}$. \refL{lem:decreasing} then implies that $\bE[I_sI_t] \leq 1/2-c$ whenever $\ip{s}{t} \geq \floor{n/3}$.

Our last ingredient is an upper bound bound on the sum 
\[
\sum_{\underset{\ip{s}{t} < \floor{n/3}}{t \in H}} w(t),
\]
for every $s \in H$.
\begin{lemma}\label{lem:small-weight}
Let $w$ be a balanced weight function.
Then for all $s \in H$ we have that 
\[
\sum_{\underset{\ip{s}{t} < \floor{n/3}}{t \in H}}w(t) \leq 3/4.
\]
\end{lemma}
The proof of this lemma is included in the appendix. 
With these lemmas we are ready to bound the variance:
\begin{theorem}\label{thm:var}
 Let $z$ be distributed uniformly over $H$. Then 
\[
\Var(\vol{\cup_{s \in S_1}K_s}) \leq 1/4 - c,
\]
for some absolute constant $c > 0$.
\end{theorem}
\begin{proof}
First note that $\bE[\vol{\cup_{s \in S_1}K_s}]^2 = (1/2 - O(n^{-1/2}))^2 = 1/4 - O(n^{-1/2})$. Next we write
\[
\bE[\vol{\cup_{s \in S_1}K_s}^2] = \sum_{s \in H}\sum_{\underset{\ip{s}{t} < \floor{n/3}}{t \in H}}w(s)w(t)\bE[I_sI_t] + \sum_{s \in H}\sum_{\underset{\ip{s}{t} \geq \floor{n/3}}{t \in H}}w(s)w(t)\bE[I_sI_t].
\]
In the first double summation we can replace $\bE[I_sI_t]$ with $1/2$ since $\bE[I_sI_t] \leq \bE[I_s] \leq 1/2$. By \refL{lem:decreasing} and \refL{lem:fluctuation} we may replace $\bE[I_sI_t]$ in the second summation by $(1/2 - c)$, with $c$ being the absolute constant from \refL{lem:fluctuation}.  This reduces the summation to 
\[
\sum_{s \in H}\sum_{\underset{\ip{s}{t} < \floor{n/3}}{t \in H}}w(s)w(t)\bE[I_sI_t] + \sum_{s \in H}\sum_{\underset{\ip{s}{t} \geq \floor{n/3}}{t \in H}}w(s)w(t)\bE[I_sI_t] \leq \sum_{s \in H} w(s) \lrpar{\sum_{\underset{\ip{s}{t} < \floor{n/3}}{t \in H}} (1/2)w(t) + \sum_{\underset{\ip{s}{t} \geq \floor{n/3}}{t \in H}} (1/2 - c)w(t)}.
\]
By \refL{lem:small-weight} the sum of weights of points in $H$ that are distance less than $\floor{n/3}$ from $s$ is at most 3/4. Since the sum of the weights of all points is 1 we may further bound the summation by
\[
\sum_{s \in H} w(s) \lrpar{\sum_{\underset{\ip{s}{t} < \floor{n/3}}{t \in H}} (1/2)w(t) + \sum_{\underset{\ip{s}{t} \geq \floor{n/3}}{t \in H}} (1/2 - c)w(t)} \leq \sum_{s \in H}w(s)\lrpar{\frac{1}{2} \cdot \frac{3}{4} + \lrpar{\frac{1}{2} -  c} \cdot \frac{1}{4}} = \frac{1}{2} - \frac{c}{8}.
\]
Therefore
\begin{align*}
\Var(\vol{\cup_{s \in S_1}K_s}) &= \bE[\vol{\cup_{s \in S_1}K_s}^2] - \bE[\vol{\cup_{s \in S_1}K_s}]^2, \\
&\leq \frac{1}{2} - \frac{c}{8} - \lrpar{\frac{1}{4} - O(n^{-1/2})}, \\ 
&= \frac{1}{4} - \frac{c}{8} + O(n^{-1/2}).
\end{align*}
and we conclude as desired.
\end{proof}
\begin{proof}[Proof of \refT{thm:uni-ub}]
For convenience we write $A := \vol{\cup_{s \in S_1}K_s},B := \vol{\cup_{s \in S_2}K_s}, C := \vol{\cup_{s \in S_3}K_s}$. We first consider the case where $n \geq n_0$ where $n_0$ is a sufficiently large constant.
By definition of $A$ and \refL{lem:sterling} $1/2 - 100/\sqrt{n} \leq \bE[A] \leq 1/2$. 
By \refT{thm:var} $\Var(A) \leq 1/4 - c$. By Chebyshev's inequality 
\[
    \bP[|A - \bE[A]| > 1/2-c] \leq \frac{1/4 - c}{1/4 - c + c^2} =: 1 - \delta.
\]
 By \refL{lem:sterling} and Markov's inequality $\bP[C > 100/(\delta\sqrt{n})] \leq \delta/2$. Therefore with probability at least $1 - (1-\delta) - \delta/2 = \delta/2$, $A,B$ and $C$ satisfy 
 \begin{align*}
 C &\leq 100/(\delta\sqrt{n}) = O(n^{-1/2}), \\
 A &\geq c - 100/\sqrt{n} = \Omega(1), \\
 B &\geq 1 - (100/(\delta\sqrt{n}) + 1 - c) = c - 100/(\delta\sqrt{n}) = \Omega(1). 
 \end{align*}
 Since the probability is non-zero there exists realizations of $S_1,S_2,S_3$ satisfying the above guarantees. In particular $\psi_K = O(n^{-1/2})$. For the case of $2 \leq n < n_0$ we use an ad-hoc argument. We proceed by casework:
 \begin{enumerate}
     \item There exists elements $s,t$ in $H$ such that $s_1 = 0, t_1 = 1$ and $w(s),w(t) > 1/4$: In this case $\ip{s}{t} \geq 2$. Otherwise (since $n \geq 2$) they agree on 1 coordinate and the sum of their weights exceeds 1/2, contradicting the balance conditions. Thus we may take $S_1 := \{s\}, S_2 := \{t\}$ in which case $\psi(K) < 2 = O(n^{-1/2})$.
     \item There does not exist an element $s$ in $H$ satisfying $s_1 = 0$ having weight greater than $1/4$: In this case let $s$ be the heaviest element with first coordinate equal to 1. By the balance condition $w(s) \geq 2^{1-n}$. Let $t,v$ be the two heaviest elements with first coordinate equal to 0. By assumption and the balanced conditions $w(t),w(v)$ are both between $2^{-n-1}$ and $1/4$. Since $t$ and $v$ are distinct and disagree with $s$ on the first coordinate without loss of generality $t$ disagrees with $s$ on at least 2 coordinates. Therefore we may take $S_1 = \{s\}, S_2 = \{t\}$ in which case $\psi_K \leq 2^{n_0+1} = O(n^{-1/2})$.
     \item There does not exist an element $s$ in $H$ satisfying $s_1 = 1$ having weight greater than $1/4$: The argument for this case is identical to the argument for the second case.
 \end{enumerate}
\end{proof} 
\section{Closing remarks and future directions}
In conclusion we have determined the $\ell_0$ isoperimetric coefficient for cubes and provided a best possible general upper-bound for the isoperimetric coefficient of any measurable set.  Surprisingly, most of the methods employed were combinatorial in nature (with some mild geometric interpretation as a source of direction). We wonder if more traditional techniques from geometric analysis would help in this line of work. \\
We leave now two open questions: 
\begin{enumerate}
\item Is it possible to apply the shaking method to other measurable sets to prove lowerbounds on the $\ell_0$ isoperimetric coefficient? We believe one such candidate of sets are the class of unconditional bodies. In particular, after applying Steiner symmetrization along all $n$ coordinate directions to axis-disjoint sets $S_1$, and $S_2$, the sets are still axis-disjoint and the distribution of $S_1$,$S_2$ within every orthant of the unconditional body will look similar to the distribution of axis-disjoint subsets after compression on the continuous unit cube. One cannot apply Harper's theorem to each orthant since it need not be a cube, but it is reasonable to conjecture that an analogous theorem should hold.

\item In \refL{lem:small-weight} we showed, with respect to a balanced weight function, that for any $s \in H$ the sum of the weights of points $t$ in $H$ satisfying $\ip{s}{t} < \floor{n/3}$ cannot be too close to 1 for. We conjecture that a similar result should hold for $\ip{s}{t} < \floor{n/2}$ if one averages over all $s \in H$. We leave the question as an open problem.
\begin{problem}
Let $w$ be a ``balanced'' function. Is it true that 
\[
\sum_{s \in H}\sum_{\underset{\ip{s}{t} < \floor{n/2}}{t \in H}}w(s)w(t) \approx \sum_{s \in H}\sum_{\underset{\ip{s}{t} > \ceil{n/2}}{t \in H}}w(s)w(t)?
\]
\end{problem}
Note that the number of terms appearing on the left hand side is equal to the number of terms appearing on the right hand side. Our proof technique from \refL{lem:small-weight} does not work here as the bound on $\sum_{\underset{\ip{s}{t}_1 < k}{t \in H}}w(t)$ grows to $1 - o(1)$ once $k = n/2 - o(n)$. However we could not come up with a ``balanced'' weight function for which the left hand side was asymptotically larger than the right hand side. 
\end{enumerate}
\section{Acknowledgements}
The author would like to thank Galyna Livshyts for bringing this problem to his attention. The author would like to thank Christopher Dupre and Galyna Livshyts for numerous helpful discussions. The author would like to thank Galyna Livshyts and Santosh Vempala for valuable feedback regarding this work. The author would like to thank Vladimir Koltchinskii for the insightful generalization of \refT{thm:cube} in the form of \refC{cor:cube}. Lastly, the author would like to thank the referees for providing valuable remarks and suggestions. This material is based upon work partially supported by a Georgia Tech ARC-ACO Fellowship and NSF-BSF DMS-2247834.
\section{Appendix}
\subsection{Proof of \refL{lem:transform-prop}}

\begin{proof}
We address the properties in order:
\begin{enumerate}
\item Our argument is an adaptation of the proof of measurability of Steiner symmetrization as appearing in \cite{measurereference}. We will use the fact that the region bounded between the graphs of two measurable functions is measurable (the proof of which can also be found in \cite{measurereference}).
Given $x \in \R^n, v \in \bS^{n-1}$ we define $\ell_v(x) := \{y~:~ y = x+tv \text{ for some }t \in \R\}$. Up to permutation of indices we may assume that $i = 1$.
By Fubini's theorem the function
 $: \R^{n-1} \rightarrow [0,\infty)$ defined by $f(b) = |A \cap \ell_v(b)|$ is measurable. Thus $\hat{A} = \{(y,b)~|~0 \leq y \leq f(b), b \in \R^{n-1}\}$ is  the region bounded between $f$ and the zero-function (both of which are measurable functions) and is therefore a measurable set. Similarly 
$\hat{B} = \{(y,b)~|~1 - f(b) < y \leq 1,b \in \R^n\} \setminus \{(1,b)~|~f(b) = 0 \}$ is the difference between a region bounded between two measurable functions (hence measurable) and a set of measure 0 (hence measurable). Thus $\hat{B}$ is measurable.
\item By Fubini-Tonelli we may write 
\begin{align*}
\vol{A} &= \int_{[0,1]^{n-1}}\int_{[0,1]} \textbf{1}_{A}(x,y)~dy~dx, \\
&= \int_{[0,1]^{n-1}} \vold{1}{y \in A~|~\Proj_{e_i^\top}y = x}~dx, \\
&= \int_{[0,1]^{n-1}} \vold{1} {[0, \vold{1} {\{ y \in A~|~\Proj_{e_i^\top}y = x \}}]}   ~dx, \\
&= \int_{[0,1]^{n-1}}\int_{[0,1]} \textbf{1}_{\hat{A}}(x,y)~dy~dx= \vol{\hat{A}}.
\end{align*}
The same argument holds for $B$.
\item Suppose $\hat{A},\hat{B}$ are not axis-disjoint. Since $A,B$ were axis-disjoint and $\Proj_{e_i\top}\hat{A} \subseteq \Proj_{e_i^{\top}}A$, $\Proj_{e_i\top}\hat{B} \subseteq \Proj_{e_i^{\top}}B$ there exists points $a \in \hat{A},b \in \hat{B}$ such that $a$ and $b$ agree on $n-1$ indices, one of which is $i$. Let $j$ denote the single index they do not agree on. Without loss of generality suppose that $i = 1, j = 2$. By the definition of $\hat{A},\hat{B}$ there exists $j_1,j_2$ such that $[0,a_1] \times a_2  \times \{a_3,\cdots,a_n\} \subseteq \hat{A}$ and $[b_1 - \eps, 1] \times b_2 \times \{b_3, \cdots, b_n\} \subseteq \hat{B}$ for some $\eps > 0$. But $a_1 = b_1$ means $b_1 - \eps < a_1$ meaning $\vold{1}{[0,a_1]} + \vold{1}{[b_1 - \eps,1]} > 1$.
Since the sum of Lebesgue measures of disjoint measurable subsets of $[0,1]$ can't exceed 1 it follows that $A \cap \{x_2 = a_2,\cdots, x_n = a_n\}$ and $B \cap \{x_2 = b_2, \cdots, x_n = b_n\}$ are not axis-disjoint, meaning $A,B$ are not axis-disjoint, a contradiction.  
\item We prove this for $\hat{A}$ (the result for $\hat{B}$ is analogous). Without loss of generality assume that $i = 1$. Since $A$ is measurable its $(n-1)$ dimensional coordinate projections are Lebesgue measurable and it immediately follows that $\cover{A}$ is Lebesgue measurable. The same is true for $A$. By Fubini-Tonelli we may write $\vol{\cover{\hat{A}}}$ and $\vol{\cover{A}}$ as integrals over the set of fibers $f$ along the $e_1$ direction. Thus it suffices to show that $\vold{1}{f \cap \cover{\hat{A}}} \leq \vold{1}{f \cap \cover{A}}$.
If $f \cap \hat{A} \neq \emptyset$ then $f \cap A \neq \emptyset$ which in particular means that $\textbf{1}_{\cover{A}} = \textbf{1}_{\cover{\hat{A}}} = 1$ over $f$. Thus assume $f \cap \hat{A} = \emptyset$. For every $j \in \{2,\cdots,n\}$ let $F_{1,j}$ denote the 2-dimensional coordinate plane that agrees with $f$ on all but the first and $j$th coordinate. Observe that any point in $\hat{A}$ that covers a portion of $f$ must be contained in some $F_{1,j}$. Furthermore the subset of $f$ that is not covered by $F_{1,j}$ is of the form $f \cap \{x_1 > s_j\}$ or $f \cap \{x_1 \geq s_j\}$ where 
\[
s_j = \sup_{\underset{w_k = f_k~\forall k \not \in \{1,j\}}{w \in \hat{A}}} w_1.
\]
Therefore the measure of the portion of the fiber that is covered by $\hat{A}$ is equal to $\max_{j \in [n]\setminus \{1\}} s_j$.
But 
\[
s_j = \sup_{y \in \R}\vold{1}{A \cap \ell_y}.
\]
where $\ell_y := \{x \in \R^n~|~ x_j = y, x_i = f_i \text{ for all } i \not \in \{1,j\}\}$. In particular the measure of the portion of the fiber covered by $\hat{A}$ is at most the measure of the portion of the fiber covered by $A$. Thus $\vold{1}{f \cap \cover{\hat{A}}} \leq \vold{1}{f \cap \cover{A}}$.
\end{enumerate}
\end{proof}
\subsection{Proof of \refL{lem:step 3}}
\begin{proof}
Recall that we may assume that $p \leq 1/2$. For convenience we write $G(k) := \bP(X=k)$. We first consider the case where $n \leq n_0$ where $n_0$ is a sufficiently large constant. Note that
\[
\frac{G(k+1)}{G(k)} = \frac{n-k}{k+1} \cdot \frac{1-p}{p} \geq \frac{1}{2np} \geq \frac{1}{n}.
\]
Therefore 
\[
\bP[X = k] \geq \frac{\Pr[X < k]}{2(n+1)^kp} \geq \frac{\min\{\bP[X < k],\bP[X > k]}{2(n_0+1)^{n_0}p}.
\]
We conclude this case by observing that $(2(n_0+1)^{n_0}p)^{-1} \geq C/\sqrt{np(1-p)}$ for some absolute constant $C$.
We now consider the case of $n > n_0$. 
Our proof strategy will depend on $p$. Suppose $p < 1/(4n)$. Then
\[
\frac{G(k+1)}{G(k)} = \frac{n-k}{k+1} \cdot \frac{1-p}{p} \geq \frac{1}{2np} > 2.
\]

Therefore for all $k$ we have 
\[
\bP[X < k] \leq 2\bP[X = k-1] \leq 4np\bP[X = k].
\]
Since $4np \leq C\sqrt{np(1-p)}$ for some absolute constant $C$ we conclude that
\[
\bP[X = k] \geq \frac{C}{\sqrt{np(1-p)}}\min\left\{\bP[X < k], \bP[X > k]\right\}.
\]
Suppose now that $1/(4n) \leq p \leq 1/2$.
 \begin{lemma}\label{lem:growth-rate}
Let $0 < x < 1$ and define $s := 1 + x$. Then 
\[
\frac{G(k+1)}{G(k)} \geq s \iff k \leq \frac{n(1-p) - sp}{(s-1)p + 1},~~\frac{G(k+1)}{G(k)} \leq \frac{1}{s} \iff k \geq \frac{sn(1-p) - p}{p + s(1-p)}.
\]
\begin{proof}
We have that 
\begin{align*}
\frac{G(k+1)}{G(k)} = \frac{\binom{n}{k+1}p^{n-{k+1}}(1-p)^{k+1}}{\binom{n}{k}p^{n-k}(1-p)^k} = \frac{n-k}{k+1} \cdot \frac{1-p}{p}.
\end{align*}
The result is now obtained by solving the corresponding inequalities in terms of $k$.
\end{proof}
\end{lemma}
Next we define $x = c_1/\sqrt{np(1-p)}$ where $c_1$ is an absolute constant which is sufficiently small. In view of $1/(4n) \leq p < 1/2$ we may assume that $x < c_2$, where $c_2$ is a sufficiently small absolute constant.
By \refL{lem:growth-rate} $\{k:~ G(k)/G(k-1) \geq 1 + x\}$ is an initial segment of $\{0,1,2,\cdots,n\}$. Taking $r$ to be any element in the set we have that 
\[
\frac{G(r)}{G(0) + \cdots + G(r-1)} \geq \frac{G(r)}{\frac{G(r)}{(1+x)^{r}} + \cdots + \frac{G(r)}{(1+x)}} \geq \frac{1}{\frac{1}{1- \frac{1}{1+x}}-1} = \frac{1}{\frac{x+1} {x}-1} = \frac{x}{x+1-x} = x.
\]
Similarly, by \refL{lem:growth-rate} $\{k:~ G(k)/G(k+1) \leq 1/(1+x)\}$ is a terminal segment of $\{0,1,2,\cdots,n\}$. Taking $r$ to be any element in the set we have that 
\[
\frac{G(r)}{G(r+1) + \cdots + G(n)} \geq \frac{G(r)}{\frac{G(r)}{(1+x)} + \cdots + \frac{G(r)}{(1+x)^{n-r-1}}} \geq \frac{1}{\frac{1}{1- \frac{1}{1+x}}-1} = \frac{1}{\frac{x+1} {x}-1} = \frac{x}{x+1-x} = x.
\]
Therefore letting $k_1 := \max\{k:~ G(k)/G(k-1) \geq 1 + x\}$ and $k_2 := \min\{k:~ G(k)/G(k+1) \leq 1/(1+x)\}$ we have that
\[
\bP[X = k] \geq x\min\left\{\bP[X < k], \bP[X > k]\right\},
\]
for all $k \not \in (k_1,k_2)$. It remains to consider $k \in (k_1,k_2)$.
For this we will bound the size of $k_2-k_1$. Observe, by our definitions of $k_1$ and $k_2$ that 
\begin{align*}
k_2 - k_1 &\leq \frac{(1+x)n(1-p) - 1}{p + (1+x)(1-p)} - \frac{n(1-p) - (1+x)p}{1 + xp}, \\ &=  \frac{((1+x)n(1-p) - 1)(1 + xp) - (n(1-p) - (1+x)p)(p + (1+x)(1-p))}{(1 + xp)(p + (1+x)(1-p))}.
\end{align*}
The denominator of this expression satisfies
\[
(1 + xp)(p + (1+x)(1-p)) \ge (1+xp)(1 + x(1-p)) \ge 1.
\]
The numerator satisfies
\begin{align*}
    ((1+x)n(1-p) - 1)(1 + xp) - (n(1-p) - (1+x)p)(p + (1+x)(1-p)), \\
    = (1+x)n(1-p)(1+xp) - n(1-p)((1+x)(1-p)+p) + \underbrace{(1+x)p(p(1+x)(1-p)) - (1+xp)}_{=:R}, \\
    = xnp(1-p)(2+x) + R, \\
    \leq \frac{2+c_2}{c_1x}.
\end{align*}
where in the last line we used the fact that $R$ can be made non-positive by taking $c_2$ sufficiently small. Altogether we conclude that 
\[
k_2-k_1 \leq \frac{2 + c_2}{c_1x} =: C_1/x,
\]
where $C_1$ is an absolute constant. 
By the definition of $k_1,k_2$ and our upperbound on $k_2-k_1$ we have that 
\[
\max_{i,j \in (k_1,k_2)} \frac{G(i)}{G(j)} \leq (1+x)^{k_2-k_1} \leq \exp((k_2-k_1)x) = \exp(C_1) =: C_2,
\]
where $C_2$ is an absolute constant. Lastly, let $k \in (k_1,k_2)$ and write $k = k_1 + t$. Then
\begin{align*}
\bP[X < k] &\leq \bP[X < k_1] + \bP[X \in [k_1,k)], \\
&\leq \frac{1}{x}\bP[X = k_1] + tC_2\bP[X = k], \\
&\leq \frac{C_2}{x}\bP[X = k] + \frac{C_1C_2}{x}\bP[X = k], \\
&= \frac{C_1C_2 + C_2}{c_1}\sqrt{np(1-p)}\bP[X=k], \\
&=: C_3\sqrt{np(1-p)}\bP[X=k].
\end{align*}
where $C_3$ is an absolute constant. Therefore when $k \in (k_1,k_2)$ we have that 
\[
\bP[X = k] \geq \frac{1}{C_3\sqrt{np(1-p)}}\bP[X < k] \geq \frac{1}{C_3\sqrt{np(1-p)}}\min\left\{\bP[X < k],\bP[X > k]\right\}.
\]
\end{proof}
\subsection{Proof of \refL{lem:fluctuation}}
\begin{proof}
For convenience we will assume that $n$ is divisible by 6 (i.e. $n = 6r$) The proof is analogous for other $n$ (albeit with messier calculations due to exactness). 
In this case the length of $z_1$ is $4r$ and the length of $z_2$ is $2r$. Recall that 
\[
\cE := \{\norm{z}_1 < \floor{n/2}\}, \cF := \{\norm{z_1}_1 - \norm{z_2}_1 < \floor{n/3} - \floor{n/3}\}.
\]
Since $z$ is distributed uniformly over $H$ and at most half of the elements in $H$ have fewer than $n/2$ ones we have that $\bP[\cE] \leq 1/2$. Next,
consider now the following two events:
\begin{align*}
    \cH_1 &:= \{ \norm{z_1}_1 \in [2r - \sqrt{r},2r]\}. \\
    \cH_2 &:= \{\norm{z_2}_1 \in [r - 3\sqrt{r}, r - 2\sqrt{r}]\}.
\end{align*}
Since $z$ is distributed uniformly on $H$ the events $\cH_1$ and $\cH_2$ are independent.
When $\cH_1 \cap \cH_2$ holds it follows that $\norm{z}_1 \leq 3r - \sqrt{r} < 3r$ and $\norm{z}_1 - \norm{z}_2 \geq r + \sqrt{r}$.
Since $\floor{n/2} = 3r$ and $\floor{n/3} = 2r$ we conclude that $\cH_1 \cap \cH_2 \subseteq \cE \cap \cF^c$. We now lowerbound $\bP[\cH_1 \cap \cH_2]$. To do so we will use the following version of the Berry-Esseen inequality for binomial random variables \cite{nagaev2011bound}:
\begin{lemma}\label{lem:BE}
Let $X$ denote a binomial random variable with parameters $n,p$ and let $\Phi$ denote the cumulative distribution function of the standard normal random variable. Then 
\[
\sup_{x \in \R}\left | \bP \lrpar{\frac{X - n(1-p)}{\sqrt{np(1-p)}} \leq x} - \Phi(x)\right| \leq \frac{0.4215(p^2+(1-p)^2)}{\sqrt{np(1-p)}}.
\]
\end{lemma}

Since $z$ is distributed uniformly over $H$, $z_1$ and $z_2$ are independent. Moreover $\norm{z_1}_1$ is distributed as a $(4r,0.5)$-binomial random variable and $\norm{z_2}_1$ is distributed as a $(2r,0.5)$. Therefore by \refL{lem:BE} we have that:
\begin{align*}
\bP(\cH_1) 
&= \bP\lrpar{2r-\sqrt{r} \leq \norm{z}_1 \leq 2r}, \\
&\geq \bP(\norm{z_1}_1 \leq 2r) - \bP(\norm{z_1}_1 \leq 2r-\sqrt{r}),\\
&= \bP\lrpar{\frac{\norm{z_1}_1 - 2r}{\sqrt{r}} \leq 0} - \bP\lrpar{\frac{\norm{z_1}_1 - 2r}{\sqrt{r}} \leq -1}, \\
&\geq \Phi(0) - \Phi(-1) - \frac{1}{2\sqrt{r}}.\\
\bP(\cH_2) 
&= \bP\lrpar{r-3\sqrt{r} \leq \norm{z_2}_1 \leq r - 2\sqrt{r}}, \\
&\geq \bP(\norm{z_2}_1 \leq r-2\sqrt{r}) - \bP(\norm{z_2}_1 \leq r-3\sqrt{r}),\\
&= \bP\lrpar{\frac{\norm{z_2}_1 - r}{\sqrt{r/2}} \leq -2\sqrt{2}} - \bP\lrpar{\frac{\norm{z_2}_1 - r}{\sqrt{r/2}} \leq -3\sqrt{2}}, \\
&\geq \Phi(-2\sqrt{2}) - \Phi(-3\sqrt{2}) - \frac{1}{\sqrt{r}}.
\end{align*} 
In addition since $\cH_1,\cH_2$ are non-empty we have that $\bP(\cH_1) \geq 2^{-4r}$ and $\bP(\cH_2) \geq 2^{-2r}$. Therefore
\begin{align*}
\bP(\cH_1) \geq c_1 &:=\min_{r \geq 1}\max\lrpar{2^{-4r},\Phi(-1) - \Phi(-2) - \frac{1}{2\sqrt{r}}}. \\
\bP(\cH_2) \geq c_2 &:=\min_{r \geq 1}\max\lrpar{2^{-2r},\Phi(-2\sqrt{2}) - \Phi(-4\sqrt{2}) - \frac{1}{\sqrt{r}}}.
\end{align*}
In particular $\bP[\cH_1 \cap \cH_2] = \bP[\cH_1]\bP[\cH_2] \geq c_1c_2$  where $c_1c_2$ is an absolute constant. Therefore
\[
\bP(\cE \cap \cF) = \bP(\cE) - \bP(\cE \cap \cF^c)
\leq \bP(\cE) - \bP(\cH_1 \cap \cH_2)
\leq \frac{1}{2} - c_1c_2.
\]
\end{proof}
\subsection{Proof of \refL{lem:small-weight}}
\begin{proof}
Fix $s \in H$. Let $\xor : H \times H \to H$ denote the xor function, where the $i$th bit of the output is 0 is the corresponding entries agree and 1 if they disagree. Since $w$ is a balanced weight function the weight function $\tilde{w}$ defined by $\tilde{w}(t) := w(t \xor s)$ is also a balanced weight function. In addition
\[
\sum_{\underset{\ip{s}{t} < \floor{n/3}}{t \in H}}w(s) = \sum_{\underset{\ip{\textbf{0}_n}{t} < \floor{n/3}}{t \in H}}\tilde{w}(s).
\]
We may therefore assume that $s = \textbf{0}_n$ and replace $w$ with $\tilde{w}$. By the second property of balanced weight functions we have that 
\[
\sum_{i = 1}^n \sum_{\underset{t_i = 0}{t \in H}}w(t) =
\sum_{i = 1}^n \sum_{\underset{t_i = 1}{t \in H}}w(t).
\]
 For every $t \in H$ the number of times $w(t)$ appears in the double summation on the left hand side is equal to the number of $0$s of $t$. Similarly, the number of times $w(t)$ appears in the double summation on the right hand side is equal to the number of $1$s of $t$. Therefore
\[
\sum_{k = 0}^n\sum_{\underset{\norm{t}_1 = k}{t \in H}}(n-k) \cdot w(t) = \sum_{k = 0}^n\sum_{\underset{\norm{t}_1 = k}{t \in H}}k \cdot w(t) \implies \sum_{k = 0}^n\sum_{\underset{\norm{t}_1 = k}{t \in H}}k \cdot w(t) = \frac{n}{2}.
\]
Suppose now that the sum of weights of points whose $\ell_1$ norm is less than $\floor{n/3}$ is larger than $3/4$. Since the sum of the weights of all points is 1 we get
\begin{align*}
\sum_{k = 0}^n\sum_{\underset{\norm{t}_1 = k}{t \in H}}k \cdot w(t) &\leq \sum_{\underset{\norm{t}_1 < \floor{n/3}}{t \in H}}(\floor{n/3}-1) \cdot w(t) + \sum_{\underset{\norm{t}_1
\geq  \floor{n/3}}{t \in H}}n \cdot w(t), \\
&< (\floor{n/3} - 1)\cdot \frac{3}{4} + n \cdot \frac{1}{4}, \\
&< \frac{n}{3} \cdot \frac{3}{4} + \frac{n}{4} = \frac{n}{2},
\end{align*}
which is a contradiction. Since $s = \textbf{0}_n$ we conclude that
\[
\sum_{\underset{\ip{s}{t} < \floor{n/3}}{t \in H}} w(t) =
\sum_{\underset{\norm{t}_1 < \floor{n/3}}{t \in H}} w(t) < 3/4. 
\]
\end{proof}
\bibliographystyle{plain}

\normalsize
\end{document}